\documentclass[letterpaper, 10pt, conference]{ieeeconf}      

\IEEEoverridecommandlockouts                              
\usepackage{graphics} 
\usepackage{comment}
\usepackage{epsfig} 
\usepackage{mathptmx} 
\usepackage{times} 
\usepackage{amsmath} 
\usepackage{amssymb}  
\usepackage{mathtools}
\usepackage{xcolor}
\usepackage{subfigure}
\usepackage{url}

\renewcommand{\baselinestretch}{0.98} 
\newtheorem{theorem}{Theorem}
\newtheorem{definition}{Definition}

\title{\LARGE \bf
Nash Equilibrium Seeking for Noncooperative Duopoly Games via Event-Triggered Control}

\author{Victor Hugo Pereira Rodrigues$^{1}$, Tiago Roux Oliveira$^{2}$, \\ Miroslav Krsti{\' c}$^{3}$, Tamer Ba{\c s}ar$^{4}$.
\thanks{$^{1}$ V. H. P. Rodrigues is with the Program of Control and Automation Engineering, Federal Center for Technological Education Celso Suckow da Fonseca (CEFET/RJ), Rio de Janeiro – RJ, Brazil. Email: {\tt\small victor.rodrigues@cefet-rj.br}}%
\thanks{$^{2}$ T. R. Oliveira is with the Dept. of Eletronics and Telecommunication Engineering, State University of Rio de Janeiro (UERJ), Rio de Janeiro -- RJ, Brazil. Email: {\tt\small tiagoroux@uerj.br}}%
       \thanks{$^{3}$ M. Krsti{\' c} is with the Dept. of Mechanical and Aerospace Engineering, University of California -- San Diego (UCSD), La Jolla -- CA, USA. Email: {\tt\small krstic@ucsd.edu}} 
         \thanks{$^{4}$ T. Ba{\c s}ar is with the Dept. of Electrical and Computer Engineering, University of Illinois Urbana-Champaign, Urbana -- IL, USA. Email: {\tt\small basar1@illinois.edu}}%
}

\begin{document}

\maketitle
\thispagestyle{empty}
\pagestyle{empty}

\begin{abstract}
This paper proposes a novel approach for locally stable convergence to Nash equilibrium in duopoly noncooperative games based on a distributed event-triggered control scheme. The proposed approach employs extremum seeking, with sinusoidal perturbation signals applied to estimate the Gradient (first derivative) of unknown quadratic payoff functions. This is the first instance of noncooperative games being tackled in a model-free fashion integrated with the event-triggered methodology. Each player evaluates independently the deviation between the corresponding current state variable and its last broadcasted value to update the player action, while they preserve control performance under limited bandwidth of the actuation paths and still guarantee stability for the closed-loop dynamics. In particular, the stability analysis is carried out using time-scaling technique, Lyapunov's direct method and averaging theory for discontinuous systems. We quantify the size of the ultimate small residual sets around the Nash equilibrium and illustrate the theoretical results numerically on an example.
\end{abstract}

\section{INTRODUCTION}

Game theory offers a conceptual framework for analyzing social interactions among competitive players, employing mathematical models to understand strategic decision-making \cite{Tirole:1991,BZ:2018}. Its application spans various domains, including engineering systems, biological behaviors, and financial markets, making it a crucial tool across diverse fields \cite{Basar:2019,Sastry:2013,BZ2:2018}. Extensive research exists on differential games, addressing both theoretical aspects and practical implementations.

Games can be broadly classified as cooperative and noncooperative \cite{Basar:1999}. Cooperative games involve players forming enforceable agreements, whereas noncooperative games focus on individual player actions and Nash equilibria \cite{N:1951}. Nash equilibrium, a key concept in noncooperative game theory, represents a state where no player can unilaterally improve its own payoff \cite{Basar:1999}.

Efforts to achieve convergence to Nash equilibrium have been ongoing for decades \cite{LB:1987,B:1987}, including studies on learning-based update schemes \cite{ZTB:2013}. Recent work has explored extremum seeking (ES) approaches for real-time computation of Nash equilibria in noncooperative games \cite{FKB:2012}, enabling stable convergence without requiring model information \cite{KW:2000}.

In contemporary networked systems, communication plays a pivotal role, impacting control performance and system behavior \cite{ZHGDDYP:2020}. Event-Triggered Control (ETC) offers a solution to mitigate traffic congestion 
by executing control tasks, non-periodically, in response to a triggering condition designed as a function of the plant's state. This strategy reduces control effort since the control update and data communication only occur when the error between the current state and the equilibrium set exceeds a value that might induce instability \cite{BH:2013}. Research in this area spans various control-estimation designs and system complexities, addressing robustness and stability concerns \cite{IB:2010a,IB:2010,APDN:2016,ZLJ:2021,CL:2019}.

Practical engineering problems, especially in networking, can benefit from game-theoretic approaches, particularly in resource allocation and optimization \cite{Basar:2019,Sastry:2013,AB:2011}. Combining game theory with event-triggered architectures presents a promising avenue for real-time optimization in networked systems, yet literature lacks exploration of ES feedback in this context \cite{Mazo_Tabuada,Wang_Lemmon,Johansson,Xinghuo_Yu}. Fortunately, we gave a positive answer to this question in our earlier work \cite{Arxiv} by introducing solutions to the problem of designing multi-variable ES algorithms based on perturbation-based (averaging-based) estimates of the model via ETC.

Addressing this gap, this paper extends the multi-input-single-output ES algorithms to  multi-input-multi-output Nash equilibrium seeking (NES) scenarios in noncooperative games, employing a distributed ES-ETC perspective. Despite the challenges posed by decentralized control in games, our approach ensures convergence to Nash equilibrium when all players utilize ES-ETC algorithms. We provide a theoretical analysis for duopoly games, leveraging time-scaling techniques, a Lyapunov function construction and averaging methods to guarantee closed-loop stability. Numerical simulations illustrate the effectiveness of our approach, showcasing its advantages over periodic sampled-data control methods \cite{KNTM:2013,HNW:2023,ZFO:2023}.


\section{Duopoly Game with Quadratic Payoffs: General Formulation}

In a duopoly game, the optimal outcomes for players P1 and P2, denoted by $y_{1}(t) \in \mathbb{R}$ and $y_{2}(t) \in \mathbb{R}$, respectively, are not solely determined by their individual actions or decision strategies, represented by $\theta_{1}(t)\in \mathbb{R}$ and $\theta_{2}(t)\in \mathbb{R}$. Instead, there exists a set of input signals $\theta^{*}=[\theta^{*}_{1},,\theta^{*}_{2}]^T \in \mathbb{R}^{2}$ where each player's strategy optimally influences the payoff functions $J_{1}$ or $J_{2}$ of the other player. Achieving this optimal equilibrium ($\theta(t)\equiv\theta^{*}$) defines the Nash equilibrium \cite{Tirole:1991}. 

In particular, we examine duopoly games where each player's payoff function takes a quadratic form. This can be expressed as a strictly concave combination of their actions:

\begin{align}
J_{1}(\theta(t)) &= \frac{H_{11}^{1}}{2}\theta_{1}^{2}(t) + \frac{H_{22}^{1}}{2}\theta_{2}^{2}(t) +H_{12}^{1}\theta_{1}(t)\theta_{2}(t) \nonumber \\
&\quad + h_{1}^{1}\theta_{1}(t) + h_{2}^{1}\theta_{2}(t) + c_{1}, \label{eq:J1} \\
J_{2}(\theta(t)) &= \frac{H_{11}^{2}}{2}\theta_{1}^{2}(t) + \frac{H_{22}^{2}}{2}\theta_{2}^{2}(t)+ H_{21}^{2}\theta_{1}(t)\theta_{2}(t) \nonumber \\
&\quad+h_{1}^{2}\theta_{1}(t) + h_{2}^{2}\theta_{2}(t) + c_{2}, \label{eq:J2}
\end{align}

Here, $J_{1}(\theta), J_{2}(\theta) : \mathbb{R}^{2} \to \mathbb{R}$ represent the payoff functions for players 1 and 2, respectively, while $\theta_{1}(t), \theta_{2}(t) \in \mathbb{R}$ denote the players' decision variables, and $H_{jk}^{i}$, $h_{j}^{i}$, $c_{i} \in \mathbb{R}$ are constants, with $H_{ii}^{i} < 0$, for all $i, j, k \in \{1, 2\}$.

For the sake completeness, let us define the mathematical representation of the Nash equilibrium $\theta^{\ast}=[\theta^{\ast}_1\,, \theta^{\ast}_2]^T$ in a two-player game as
\begin{align} \label{Nashcu}
J_1(\theta_1^{\ast}\,,\theta_{2}^{\ast}) \geq J_1(\theta_1\,,\theta_{2}^{\ast}) ~~  \mbox{and} ~~~ J_2(\theta_1^{\ast}\,,\theta_{2}^{\ast}) \geq J_2(\theta_1^{\ast}\,,\theta_{2}). \!\!
\end{align}

Thus, no player has an incentive to independently deviate from $\theta^*$. In the duopoly context, $\theta_1$ and $\theta_2$ are elements of $\mathbb{R}$, representing the set of real numbers.

In order to determine the Nash equilibrium solution in strictly concave quadratic games involving two players, where each action set is the entire real line, one should differentiate $J_{i}$ with respect to $\theta_{i}(t) \,, \forall i=1,2$, setting the resulting expressions equal to zero, and solving the set of equations thus obtained. This set of equations, which also provides a sufficient condition due to the strict concavity, is
\begin{align}
\begin{cases}
H_{11}^{1}\theta_{1}^{*}+H_{12}^{1}\theta_{2}^{*}+h_{1}^{1}=0 \\ 
H_{21}^{2}\theta_{1}^{*}+H_{22}^{2}\theta_{2}^{*}+h_{2}^{2}=0 \\ 
\end{cases}\,,
\label{eq:NE_v0}
\end{align} 
which can be written in the matrix form as 
\begin{align}
\begin{bmatrix}
 H_{11}^{1} & H_{12}^{1} \\
 H_{21}^{2} & H_{22}^{2}
\end{bmatrix}
\begin{bmatrix}
\theta_{1}^{*} \\
\theta_{2}^{*} 
\end{bmatrix}
=-
\begin{bmatrix}
h_{1}^{1} \\
h_{2}^{2}  
\end{bmatrix}
\,.
\end{align}
Defining the Hessian matrix $H$, and vectors $\theta^*$ and $h$ by 
\begin{align}
H:=
\begin{bmatrix}
 H_{11}^{1} &  H_{12}^{1} \\
 H_{21}^{2} &  H_{22}^{2}
\end{bmatrix}
\,, \quad 
\theta^{*}:=
\begin{bmatrix}
\theta_{1}^{*} \\
\theta_{2}^{*} 
\end{bmatrix}
\,, \quad
h:=
\begin{bmatrix}
h_{1}^{1} \\
h_{2}^{2}  
\end{bmatrix}
\,, \label{eq:Htheta*h}
\end{align}
there exists an unique Nash Equilibrium at  $\theta^{*}=-H^{-1}h$, if $H$ is invertible:
\begin{align}
\begin{bmatrix}
\theta_{1}^{*} \\
\theta_{2}^{*} 
\end{bmatrix}
=-
\begin{bmatrix}
 H_{11}^{1} & H_{12}^{1} \\
 H_{21}^{2} &  H_{22}^{2}
\end{bmatrix}^{-1}
\begin{bmatrix}
h_{1}^{1} \\
h_{2}^{2}  
\end{bmatrix} \label{eq:NE_v2} \,.
\end{align}
For more details, see \cite[Chapter 4]{Basar:1999}.

\subsection{Continuous-time Extremum Seeking}

The primary goal of this paper is to develop distributed event-triggered control policies to obtain NES in noncooperative duopoly games. Throughout our analysis,
it is assumed that $H_{ii}^{i}<0$ for both $i$ (a natural consequence of strict concavity), the matrix $H$ in (\ref{eq:Htheta*h}) has full rank, and all the parameters are unknown.

\begin{figure}[ht]
\begin{center}
\includegraphics[width=3.35 in]{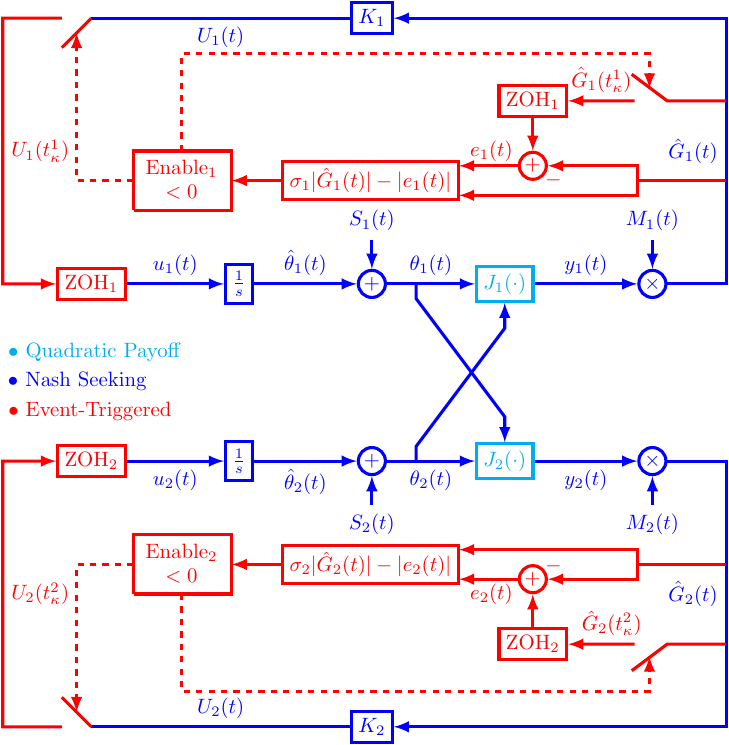}
\end{center}
\caption{Block diagram illustrating the NES strategy through distributed event-triggered control policies performed for each player.}
\label{fig:blockDiagram_v3}
\end{figure}
The schematic diagram depicted in Fig.~\ref{fig:blockDiagram_v3} provides an overview of the proposed NES policy for each player, detailing their respective outputs as follows:
\begin{align}
\begin{cases}
y_{1}(t)&=J_{1}(\theta(t))\,, \\ 
y_{2}(t)&=J_{2}(\theta(t))\,.
\end{cases}
\label{eq:yi_v0}
\end{align}
The probing signals are
\begin{align}
\begin{cases}
S_{1}(t)&=a_{1}\sin(\omega_{1}t) \\
S_{2}(t)&=a_{2}\sin(\omega_{2}t)
\end{cases}
\,,\label{eq:Si}
\end{align}
and the demodulating signals are represented by
\begin{align}
\begin{cases}
M_{1}(t)&=\frac{2}{a_{1}}\sin(\omega_{1}t) \\
M_{2}(t)&=\frac{2}{a_{2}}\sin(\omega_{2}t)
\end{cases}
\,, \label{eq:Mi}
\end{align}
having non-zero constant amplitudes $a_{1}$ and $a_{2}$, both greater than zero, and frequencies $\omega_1$ and $\omega_2$, satisfying $\omega_1 \neq \omega_2$. These probing frequencies $\omega_i$ can be chosen as
\begin{equation}\label{omegadefinition}
\omega_i=\omega_{i}'\omega=\mathcal{O}(\omega)\,, \quad
i=1 \mbox{ or }2\,,
\end{equation}
where $\omega$ represents a positive constant and $\omega_{i}'$ is a rational number. One potential selection is detailed in \cite{GKN:2012}.

If we denote $\hat{\theta}_{1}(t)$ and $\hat{\theta}_{2}(t)$ as the estimates of $\theta^{\ast}_{1}$ and $\theta^{\ast}_{2}$, respectively, we can define the ``estimation error'' as follows:
\begin{align}
\begin{cases}
\tilde{\theta}_{1}(t)=\hat{\theta}_{1}(t)-\theta_{1}^{*}\\
\tilde{\theta}_{2}(t)=\hat{\theta}_{2}(t)-\theta_{2}^{*}
\end{cases}\,.\label{eq:tildeThetai}
\end{align}
The estimation of the gradients of the unknown payoff functions 
\begin{align}
\begin{cases}
\hat{G}_{1}(t)&=M_1(t)y_1(t) \\
\hat{G}_{2}(t)&=M_2(t)y_2(t)
\end{cases}
\,, \label{eq:Gi_gold}
\end{align}
can be rewritten as
\begin{small}
\begin{align}
\hat{G}_{1}(t)&= H^{1}_{11}\theta^{\ast}_{1} +H^{1}_{12}\theta^{\ast}_{2} + h^{1}_{1} \nonumber \\
&\quad +\mathcal{H}_{11}(t)\tilde{\theta}_{1}(t) +\mathcal{H}_{12}(t)\tilde{\theta}_{2}(t) \nonumber \\
&\quad+\frac{H^{1}_{11}}{a_{1}}\sin(\omega_{1}t)\tilde{\theta}_{1}^2(t)+ \frac{H^{1}_{22}}{a_{1}}\sin(\omega_{1}t)\tilde{\theta}_{2}^2(t)  \nonumber \\
&\quad+ \frac{2H^{1}_{12}}{a_{1}}\sin(\omega_{1}t)\tilde{\theta}_{1}(t)\tilde{\theta}_{2}(t) +\delta_{1}(t)\,, \label{eq:hatG1_20240318} \\
\hat{G}_{2}(t)&= H^{2}_{21}\theta^{\ast}_{1} +H^{2}_{22}\theta^{\ast}_{2} + h^{2}_{2} \nonumber \\
&\quad +\mathcal{H}_{21}(t)\tilde{\theta}_{1}(t) +\mathcal{H}_{22}(t)\tilde{\theta}_{2}(t) \nonumber \\
&\quad+\frac{H^{2}_{11}}{a_{2}}\sin(\omega_{2}t)\tilde{\theta}_{1}^2(t)+ \frac{H^{2}_{22}}{a_{2}}\sin(\omega_{2}t)\tilde{\theta}_{2}^2(t)  \nonumber \\
&\quad+ \frac{2H^{2}_{21}}{a_{2}}\sin(\omega_{2}t)\tilde{\theta}_{1}(t)\tilde{\theta}_{2}(t) +\delta_{2}(t)\,, \label{eq:hatG2_20240318}
\end{align}
\end{small}
with time-varying parameters given by
\begin{small}
\begin{align}
\mathcal{H}_{11}(t)&=H^{1}_{11}-H^{1}_{11}\cos(2\omega_{1}t)   \nonumber \\
&\quad+\frac{2h^{1}_{1}}{a_{1}}\sin(\omega_{1}t)+ \frac{2H^{1}_{11}\theta^{\ast}_{1}}{a_{1}}\sin(\omega_{1}t) \nonumber \\
&\quad +\frac{H^{1}_{12}a_{2}}{a_{1}}\cos[(\omega_{1}+\omega_{2})t] \nonumber \\
&\quad +\frac{2H^{1}_{12}\theta^{\ast}_{2}}{a_{1}}\sin(\omega_{1}t)+\frac{H^{1}_{12}a_{2}}{a_{1}}\cos[(\omega_{1}-\omega_{2})t] \,, \\
\mathcal{H}_{12}(t)&= H^{1}_{12}- H^{1}_{12}\cos(2\omega_{1}t) \nonumber \\
&\quad+ \frac{2h^{1}_{2}}{a_{1}}\sin(\omega_{1}t)+ \frac{2H^{1}_{12}\theta^{\ast}_{1}}{a_{1}}\sin(\omega_{1}t) \nonumber \\
&\quad + \frac{2H^{1}_{22}\theta^{\ast}_{2}}{a_{1}}\sin(\omega_{1}t)+\frac{H^{1}_{22}a_{2}}{a_{1}}\cos[(\omega_{1}-\omega_{2})t] \nonumber \\
&\quad+\frac{H^{1}_{22}a_{2}}{a_{1}}\cos[(\omega_{1}+\omega_{2})t]\,, \\
\mathcal{H}_{21}(t)&= H^{2}_{21}- H^{2}_{21}\cos(2\omega_{2}t) \nonumber \\
&\quad+ \frac{2h^{2}_{1}}{a_{2}}\sin(\omega_{2}t)+ \frac{2H^{2}_{21}\theta^{\ast}_{1}}{a_{2}}\sin(\omega_{2}t) \nonumber \\
&\quad + \frac{2H^{2}_{22}\theta^{\ast}_{2}}{a_{2}}\sin(\omega_{2}t)+\frac{H^{2}_{22}a_{1}}{a_{2}}\cos[(\omega_{2}-\omega_{1})t] \nonumber \\
&\quad+\frac{H^{2}_{22}a_{1}}{a_{2}}\cos[(\omega_{2}+\omega_{1})t]\,, \\
\mathcal{H}_{22}(t)&=H^{2}_{22}-H^{2}_{22}\cos(2\omega_{2}t)   \nonumber \\
&\quad+\frac{2h^{2}_{2}}{a_{2}}\sin(\omega_{2}t)+ \frac{2H^{2}_{22}\theta^{\ast}_{2}}{a_{2}}\sin(\omega_{2}t) \nonumber \\
&\quad +\frac{H^{2}_{21}a_{1}}{a_{2}}\cos[(\omega_{2}+\omega_{1})t] \nonumber \\
&\quad +\frac{2H^{2}_{21}\theta^{\ast}_{1}}{a_{2}}\sin(\omega_{2}t)+\frac{H^{2}_{21}a_{1}}{a_{2}}\cos[(\omega_{2}-\omega_{1})t] \,,
\end{align}
\end{small}
\vspace{-1cm}
\begin{small}
\begin{align}
\delta_{1}(t)&= \frac{H^{1}_{11}a_{1}}{2}\sin(\omega_{1}t)+\frac{H^{1}_{11}a_{1}}{4}\cos(\omega_{1}t)+\frac{H^{1}_{11}a_{1}}{4}\cos(3\omega_{1}t) \nonumber \\ 
&\quad +\frac{H^{1}_{22}a_{2}^2}{2a_{1}}\sin(\omega_{1}t)-\frac{H^{1}_{22}a_{2}^2}{2a_{1}}\cos[(\omega_{1}-2\omega_{2})t] \nonumber \\
&\quad+\frac{H^{1}_{22}a_{2}^2}{2a_{1}}\cos[(\omega_{1}+2\omega_{2})t]+\frac{2H^{1}_{12}\theta^{\ast}_{1}\theta^{\ast}_{2}}{a_{1}}\sin(\omega_{1}t)   \nonumber \\
&\quad-h^{1}_{1}\cos(2\omega_{1}t)  \nonumber \\
&\quad+ \frac{a_{2}h^{1}_{2}}{a_{1}}\cos[(\omega_{1}-\omega_{2})t]+ \frac{a_{2}h^{1}_{2}}{a_{1}}\cos[(\omega_{1}+\omega_{2})t] \nonumber \\
&\quad-H^{1}_{11}\theta^{\ast}_{1}\cos(2\omega_{1}t)  \nonumber \\
&\quad + \frac{H^{1}_{12}a_{2}\theta^{\ast}_{1}}{a_{1}}\cos[(\omega_{1}-\omega_{2})t]+ \frac{H^{1}_{12}a_{2}\theta^{\ast}_{1}}{a_{1}}\cos[(\omega_{1}+\omega_{2})t]\nonumber \\
&\quad -H^{1}_{12}\theta^{\ast}_{2}\cos(2\omega_{1}t)  \nonumber \\
&\quad + \frac{H^{1}_{22}a_{2}\theta^{\ast}_{2}}{a_{1}}\cos[(\omega_{1}-\omega_{2})t] - \frac{H^{1}_{22}a_{2}\theta^{\ast}_{2}}{a_{1}}\cos[(\omega_{1}+\omega_{2})t]\nonumber \\
&\quad+\frac{2c_{1}}{a_{1}}\sin(\omega_{1}t) + \frac{2h^{1}_{1}\theta^{\ast}_{1}}{a_{1}}\sin(\omega_{1}t) + \frac{2h^{1}_{2}\theta^{\ast}_{2}}{a_{1}}\sin(\omega_{1}t)  \nonumber \\
&\quad+\frac{H^{1}_{11}\theta^{\ast 2}_{1}}{a_{1}}\sin(\omega_{1}t) + \frac{H^{1}_{22}\theta^{\ast 2}_{2}}{a_{1}}\sin(\omega_{1}t) \nonumber \\
&\quad +\frac{H^{1}_{12}a_{2}}{2}\cos(\omega_{2}t)-\frac{H^{1}_{12}a_{2}}{2}\cos[(2\omega_{1}-\omega_{2})t] \nonumber \\
&\quad+\frac{H^{1}_{12}a_{2}}{2}\cos(\omega_{2}t)+\frac{H^{1}_{12}a_{2}}{2}\cos[(2\omega_{1}+\omega_{2})t]\,,
\end{align}
\end{small}

\begin{small}
\begin{align}
\delta_{2}(t)&= \frac{H^{2}_{11}a_{2}}{2}\sin(\omega_{2}t)+\frac{H^{2}_{11}a_{2}}{4}\cos(\omega_{2}t)+\frac{H^{2}_{11}a_{2}}{4}\cos(3\omega_{2}t) \nonumber \\ 
&\quad +\frac{H^{2}_{22}a_{2}}{2}\sin(\omega_{2}t)-\frac{H^{2}_{22}a_{2}}{2}\cos[(\omega_{2}-2\omega_{1})t] \nonumber \\
&\quad+\frac{H^{2}_{22}a_{2}}{2}\cos[(\omega_{2}+2\omega_{1})t]+\frac{2H^{2}_{21}\theta^{\ast}_{1}\theta^{\ast}_{2}}{a_{2}}\sin(\omega_{2}t)   \nonumber \\
&\quad-h^{2}_{2}\cos(2\omega_{2}t)  \nonumber \\
&\quad+ h^{1}_{2}\cos[(\omega_{2}-\omega_{1})t]+ \frac{h^{2}_{2}}{a_{2}}\cos[(\omega_{1}+\omega_{2})t] \nonumber \\
&\quad-H^{2}_{11}\theta^{\ast}_{1}\cos(2\omega_{2}t)  \nonumber \\
&\quad + \frac{H^{2}_{21}a_{2}\theta^{\ast}_{1}}{a_{1}}\cos[(\omega_{2}-\omega_{1})t]+ H^{2}_{21}\theta^{\ast}_{1}\cos[(\omega_{1}+\omega_{2})t]\nonumber \\
&\quad -H^{2}_{21}\theta^{\ast}_{2}\cos(2\omega_{2}t)  \nonumber \\
&\quad + \frac{H^{2}_{22}a_{1}\theta^{\ast}_{2}}{a_{2}}\cos[(\omega_{2}-\omega_{1})t] - \frac{H^{2}_{22}a_{1}\theta^{\ast}_{2}}{a_{2}}\cos[(\omega_{2}+\omega_{1})t]\nonumber \\
&\quad+\frac{2c_{2}}{a_{2}}\sin(\omega_{2}t) + \frac{2h^{2}_{1}\theta^{\ast}_{1}}{a_{2}}\sin(\omega_{2}t) + \frac{2h^{2}_{2}\theta^{\ast}_{2}}{a_{2}}\sin(\omega_{2}t)  \nonumber \\
&\quad+\frac{H^{2}_{11}\theta^{\ast 2}_{1}}{a_{2}}\sin(\omega_{2}t) + \frac{H^{2}_{22}\theta^{\ast 2}_{2}}{a_{2}}\sin(\omega_{2}t) \nonumber \\
&\quad +\frac{H^{2}_{21}a_{1}}{2}\cos(\omega_{1}t)-\frac{H^{2}_{21}a_{1}}{2}\cos[(2\omega_{2}-\omega_{1})t] \nonumber \\
&\quad+\frac{H^{2}_{21}a_{1}}{2}\cos(\omega_{1}t)+\frac{H^{2}_{21}a_{1}}{2}\cos[(2\omega_{2}+\omega_{1})t]\,.
\end{align}
\end{small}
It is worth noting that the quadratic terms in (\ref{eq:hatG1_20240318}) and (\ref{eq:hatG2_20240318})  can be disregarded in a local analysis \cite{K:2014}. Moreover, from (\ref{eq:NE_v0}), $H_{11}^{1}\theta_{1}^{*}+H_{12}^{1}\theta_{2}^{*}+h_{1}^{1}=0$ and $H_{21}^{2}\theta_{1}^{*}+H_{22}^{2}\theta_{2}^{*}+h_{2}^{2}=0$. Therefore, the gradient estimate is locally expressed as
\begin{align}
\hat{G}_{1}(t)&= \mathcal{H}_{11}(t)\tilde{\theta}_{1}(t) +\mathcal{H}_{12}(t)\tilde{\theta}_{2}(t) +\delta_{1}(t)\,, \label{eq:hatG1_v2_20240318} \\
\hat{G}_{2}(t)&= \mathcal{H}_{21}(t)\tilde{\theta}_{1}(t) +\mathcal{H}_{22}(t)\tilde{\theta}_{2}(t) +\delta_{2}(t)\,. \label{eq:hatG2_v2_20240318}
\end{align}
\newpage

Now, considering the time-derivative of (\ref{eq:tildeThetai}) and referring to the NES scheme illustrated in Fig.~\ref{fig:blockDiagram_v3}, we can derive the dynamics governing $\hat{\theta}(t)$ and $\tilde{\theta}(t)$ as follows:
\begin{align}
\frac{d\tilde{\theta}_{1}(t)}{dt}&=\frac{d\hat{\theta}_{1}(t)}{dt}=u_{1}(t) \label{eq:dotTildeTheta_1}\,, \\
\frac{d\tilde{\theta}_{2}(t)}{dt}&=\frac{d\hat{\theta}_{2}(t)}{dt}=u_{2}(t)\,, \label{eq:dotTildeTheta_2}
\end{align}
where $u_{1}(t)$ and $u_{2}(t)$ represent the NES control laws to be formulated.

The continuous-time feedback law is given by:
\begin{align}
\begin{cases}
u_{1}(t)&=K_{1}\hat{G}_{1}(t) \\
u_{2}(t)&=K_{2}\hat{G}_{2}(t) 
\end{cases}\,,
\quad \text{for all } t\geq 0\,. \label{eq:U_continuous}
\end{align}
This controller is stabilizing, with the gain represented as:
\begin{align}
K=\begin{bmatrix} K_{1} & 0 \\
                    0   & K_{2} \end{bmatrix}\,,
\end{align}
ensuring that $KH$ is Hurwitz.

Here, control updates are triggered only for a given  sequence of time instants $(t_{\kappa})_{\kappa\in\mathbb{N}}$ determined by an event-generator. This generator is designed to maintain stability and robustness. The control task is orchestrated by a monitoring mechanism that triggers updates when the difference between the current output value and its previously computed value at time $t_\kappa$ exceeds a predefined threshold, determined by a constructed triggering condition \cite{HJT:2012}. It is important to note that in conventional sampled-data implementations, execution times are evenly spaced in time, with $t_{\kappa+1}=t_{\kappa}+ h$, where $h>0$ is a known constant, for all $\kappa\in \mathbb{N}$. However, in an event-triggered scheme, sampling times may occur aperiodically.

\subsection{Emulation of Continuous-Time Extremum Seeking Design}

We define the players' actions as:
\begin{align}
u_1(t)&=K_{1}\hat{G}_{1}(t^{1}_{\kappa}) , \quad \forall t \in \lbrack t^{1}_{\kappa}, t^{1}_{\kappa+1}) , \label{eq:u1_v1} \\
u_2(t)&=K_{2}\hat{G}_{2}(t^{2}_{\kappa}) , \quad \forall t \in \lbrack t^{2}_{\kappa}, t^{2}_{\kappa+1}) . \label{eq:u2_v1}
\end{align}
The following functions $e_{1}$ and $e_{2}$, both mapping from $\mathbb{R}$ to $\mathbb{R}$, represent the difference between the current state variable and its previously transmitted value for players P1 and P2, respectively, and are expressed as
\begin{align}
e_{1}(t)&:=\hat{G}_{1}(t^{1}_{\kappa})-\hat{G}_{1}(t) , \quad \forall t \in \lbrack t^{1}_{\kappa}, t^{1}_{\kappa+1}) , \quad \kappa\in \mathbb{N} , \label{eq:e1} \\
e_{2}(t)&:=\hat{G}_{2}(t^{2}_{\kappa})-\hat{G}_{2}(t) , \quad \forall t \in \lbrack t^{2}_{\kappa}, t^{2}_{\kappa+1}) , \quad \kappa\in \mathbb{N} . \label{eq:e2}
\end{align}

Using equations (\ref{eq:u1_v1})--(\ref{eq:e2}), we rewrite the distributed event-triggered control policies for each player in terms of $e_{1}$ and $e_{2}$, as follows
\begin{align}
u_1(t)&=K_{1}\hat{G}_{1}(t)+K_{1}e_{1}(t) , \quad \forall t \in \lbrack t^{1}_{\kappa}, t^{1}_{\kappa+1}) , \label{eq:u1_v2}\\
u_2(t)&=K_{2}\hat{G}_{2}(t)+K_{2}e_{2}(t) , \quad \forall t \in \lbrack t^{2}_{\kappa}, t^{2}_{\kappa+1}) . \label{eq:u2_v2}
\end{align}

Consequently, plugging the event-triggered control law (\ref{eq:u1_v2})--(\ref{eq:u2_v2}) into the time-derivative of (\ref{eq:hatG1_v2_20240318}), (\ref{eq:hatG2_v2_20240318}), (\ref{eq:dotTildeTheta_1}) and equation (\ref{eq:dotTildeTheta_2}), for all $t\in[t_{\kappa},t_{\kappa+1})$, where $t_{\kappa}=\min\{t^{1}_{\kappa},t^{2}_{\kappa}\}$ and $t_{\kappa+1}=\min\{t^{1}_{\kappa+1},t^{2}_{\kappa+1}\}$ (for more details, see \cite{TC:2014}), yields the following dynamics governing $\hat{G}(t)$ and $\tilde{\theta}(t)$:
\begin{align}
\frac{d\hat{G}_{1}(t)}{dt}
&= \mathcal{H}_{11}(t)K_{1}\hat{G}_{1}(t) +\mathcal{H}_{11}(t)K_{1}e_{1}(t)  \nonumber \\
&\quad+\mathcal{H}_{12}(t)K_{2}\hat{G}_{2}(t)+\mathcal{H}_{12}(t)K_{2}e_{2}(t)\nonumber \\
&\quad+\frac{d\mathcal{H}_{11}(t)}{dt}\tilde{\theta}_{1}(t) +\frac{d\mathcal{H}_{12}(t)}{dt}\tilde{\theta}_{2}(t)+\frac{d\delta_{1}(t)}{dt}\,, \label{eq:dotHatG1_v1_20240319} \\
\frac{d\hat{G}_{2}(t)}{dt}
&= \mathcal{H}_{21}(t)K_{1}\hat{G}_{1}(t) +\mathcal{H}_{21}(t)K_{1}e_{1}(t)  \nonumber \\
&\quad+\mathcal{H}_{22}(t)K_{2}\hat{G}_{2}(t)+\mathcal{H}_{22}(t)K_{2}e_{2}(t)\nonumber \\
&\quad+\frac{d\mathcal{H}_{21}(t)}{dt}\tilde{\theta}_{1}(t) +\frac{d\mathcal{H}_{22}(t)}{dt}\tilde{\theta}_{2}(t)+\frac{d\delta_{2}(t)}{dt}\,, \label{eq:dotHatG2_v1_20240319} \\
\frac{d\tilde{\theta}_{1}(t)}{dt}
&= K_{1}\mathcal{H}_{11}(t)\tilde{\theta}_{1}(t) +K_{1}\mathcal{H}_{12}(t)\tilde{\theta}_{2}(t)\nonumber \\
&\quad+K_{1}e_{1}(t) +K_{1}\delta_{1}(t)\,, \label{eq:dotTildeTheta1_v1_20240319} \\
\frac{d\tilde{\theta}_{2}(t)}{dt}
&= K_{2}\mathcal{H}_{21}(t)\tilde{\theta}_{1}(t) +K_{2}\mathcal{H}_{22}(t)\tilde{\theta}_{2}(t)\nonumber \\
&\quad+K_{2}e_{2}(t) +K_{2}\delta_{2}(t)\,. \label{eq:dotTildeTheta2_v1_20240319}
\end{align}
Hence, in a concise vector-form representation, we have
\begin{align}
\frac{d\hat{G}(t)}{dt}&= \mathcal{H}(t)K\hat{G}(t) +\mathcal{H}(t)Ke(t)+\frac{d\mathcal{H}(t)}{dt}\tilde{\theta}(t)+\frac{d\delta(t)}{dt}\,, \label{eq:dotHatG_v1_20240319} \\
\frac{d\tilde{\theta}(t)}{dt}&= K\mathcal{H}(t) \tilde{\theta}(t) +Ke(t) +K\delta(t)\,, \label{eq:dotTildeTheta_v1_20240319} 
\end{align}
where 
\begin{align}
\mathcal{H}(t)=\begin{bmatrix} 
									\mathcal{H}_{11}(t) & \mathcal{H}_{12}(t) \\
									\mathcal{H}_{21}(t) & \mathcal{H}_{22}(t)
							 \end{bmatrix} \mbox{ and } 
\delta(t)=\begin{bmatrix} 
									\delta_{1}(t) \\
									\delta_{2}(t)
							 \end{bmatrix} \,.
\end{align}
The closed-loop system described by (\ref{eq:dotHatG_v1_20240319}) and (\ref{eq:dotTildeTheta_v1_20240319}) highlights a crucial point: while the product $\mathcal{H}(t)K$ on averaging sense forms a Hurwitz matrix, the convergence to the equilibrium $\hat{G}\equiv0$ and $\tilde{\theta}\equiv0$ is not guaranteed due to the presence of the error vector $e(t):=[e_{1}(t)\,,e_{2}(t)]^{T}$ and the time-varying term $\delta(t)$ and their derivatives. However, the system does exhibit Input-to-State Stability (ISS) concerning the error vector $e(t)$ and such time-varying disturbances. Additionally, it is important to note that the disturbances $\delta(t)$ and $\frac{d\delta(t)}{dt}$ as well as the time-varying matrix $\frac{d\mathcal{H}(t)}{dt}$ possess zero mean values.

In the next two sections, we introduce the static triggering mechanism for NES, as outlined in Definitions~\ref{def:staticEvent} and \ref{def:averageStaticEvent}. This mechanism represents a fusion of distributed event-triggered data transmission with an extremum-seeking control system.

\section{Distributed Event-Triggered Actions for Nash Equilibrium Seeking}

Definition~\ref{def:staticEvent} in sequence outlines the utilization of small parameters $\sigma_{i}$, along with the errors $e_{i}$ representing the disparity between a player's current state variable and its last broadcasted value, and measurements of the gradient estimate $\hat{G}_{i}$. These components are employed to construct the ``NES Static Triggering Condition''. This approach involves updating the control laws (\ref{eq:u1_v1}) and (\ref{eq:u2_v1}) and the ZOH actuators, as depicted in Fig.~\ref{fig:blockDiagram_v3}. This ensures the asymptotic stability of the closed-loop system.

\begin{definition}[\small{NES Static-Triggering Condition}] \label{def:staticEvent}
The NES-based event-triggered controller with static-triggering condition consists of two components:
\begin{enumerate}
	\item The set of increasing sequences of time $I=\{I_{1}\,,I_{2}\}$ such that $I_{i}=\{t^{i}_{0}\,, t^{i}_{1}\,, t^{i}_{2}\,,\ldots\}$ with $t^{i}_{0}=0$, for all $i\in \{1,2\}$, generated under the following rules:
		\begin{itemize}
			\item If $\left\{t \in\mathbb{R}^{+}: t>t^{i}_{\kappa} ~ \wedge ~ \sigma_{i}|\hat{G}_{i}(t)| -|e_{i}(t)| < 0 \right\} = \emptyset$, then the set of the times of the events is $I_{i}=\{t^{i}_{0}\,, t^{i}_{1}\,, \ldots, t^{i}_{\kappa}\}$.
			\item If $\left\{t \in\mathbb{R}^{+}: t>t^{i}_{\kappa} ~ \wedge ~ \sigma_{i}|\hat{G}_{i}(t)| -|e_{i}(t)| < 0 \right\} \neq \emptyset$, the next event time is given by
				\begin{align}
					 t^{i}_{\kappa+1}\!&=\!\inf\left\{t \in\mathbb{R}^{+}: t\!>\!t^{i}_{\kappa} \wedge \sigma_{i}|\hat{G}_{i}(t)| -|e_{i}(t)| \!<\! 0 \right\}, \label{eq:tk+1_event}
				\end{align}
				\!\!\! consisting of the static event-triggering mechanism.
		\end{itemize}
	\item The $i$-th feedback control action updated at the  triggering instants in (\ref{eq:u1_v1}) and (\ref{eq:u2_v1}).
\end{enumerate}  
\end{definition}

\section{Closed-loop System for Time-Scaled Triggering Mechanism} \label{sec:closed_loop}

\subsection{Time-Scaling Procedure} \label{sec:time_scaling}

To facilitate the stability analysis of the closed-loop system, we introduce a convenient time scale. By examining (\ref{omegadefinition}), it is evident that the dither frequencies (\ref{eq:Si}) and (\ref{eq:Mi}), as well as their combinations, are rational. Additionally, there exists a time period $T$ such that
\begin{align}
T&= 2\pi \times \text{LCM}\left\{\frac{1}{\omega_{1}},\frac{1}{\omega_{2}}\right\}\,,\label{eq:period_T}
\end{align}
where LCM denotes the least common multiple. Hence, we define a new time scale for the dynamics (\ref{eq:dotHatG_v1_20240319}) and (\ref{eq:dotTildeTheta_v1_20240319}) using the transformation $\bar{t}=\omega t$, where
\begin{align}
\omega&:=\frac{2\pi}{T}\,. \label{eq:new_omega}
\end{align}
Consequently, the system (\ref{eq:dotHatG_v1_20240319}) and (\ref{eq:dotTildeTheta_v1_20240319}) can be expressed as, for all $t \in \lbrack t_{\kappa},, t_{\kappa+1})$, and $\kappa\in \mathbb{N}$:
\begin{align}
\frac{d\hat{G}(\bar{t})}{d\bar{t}}&= \frac{1}{\omega}\mathcal{H}(\bar{t})K\hat{G}(\bar{t}) +\frac{1}{\omega}\mathcal{H}(\bar{t})Ke(\bar{t})+\nonumber \\
&\quad+\frac{d\mathcal{H}(\bar{t})}{d\bar{t}}\tilde{\theta}(\bar{t})+\frac{1}{\omega}\frac{d\delta(\bar{t})}{d\bar{t}}\,, \label{eq:dotHatG_v2_20240319} \\
\frac{d\tilde{\theta}(\bar{t})}{d\bar{t}}&= \frac{1}{\omega}K\mathcal{H}(\bar{t})\tilde{\theta}(\bar{t}) +\frac{1}{\omega}Ke(\bar{t}) +\frac{1}{\omega}K\delta(\bar{t})\,. \label{eq:dotTildeTheta_v2_20240319} 
\end{align}
Now, we can implement a suitable averaging mechanism within the transformed time scale $\bar{t}$ based on the dynamics (\ref{eq:dotHatG_v2_20240319}) and (\ref{eq:dotTildeTheta_v2_20240319}). Despite the non-periodicity of the triggering events and discontinuity on the right-hand side of equations (\ref{eq:dotHatG_v2_20240319}) and (\ref{eq:dotTildeTheta_v2_20240319}), the closed-loop system maintains its periodicity over time due to the periodic probing and demodulation signals. This unique characteristic allows for the application of the averaging results established by Plotnikov \cite{P:1979} to this particular setup. 

\subsection{Average Closed-Loop System}

Let us define the augmented state as:
\begin{align}
X^{T}(\bar{t}):=\begin{bmatrix} \hat{G}(\bar{t})\ \tilde{\theta}(\bar{t}) \end{bmatrix}\,,
\end{align}
which reduces the system (\ref{eq:dotHatG_v2_20240319}) and (\ref{eq:dotTildeTheta_v2_20240319}) to:
\begin{align}
\dfrac{dX(\bar{t})}{d\bar{t}}&=\dfrac{1}{\omega}\mathcal{F}\left(\bar{t},X,\dfrac{1}{\omega}\right)\,. \label{eq:dotX_event}
\end{align}

The system (\ref{eq:dotX_event}) features a small parameter $1/\omega$ and a $T$-periodic function $\mathcal{F}\left(\bar{t},X,\dfrac{1}{\omega}\right)$ in $\bar{t}$. Therefore, the averaging theorem \cite{P:1979} can be applied to $\mathcal{F}\left(\bar{t},X,\dfrac{1}{\omega}\right)$ at $\displaystyle \lim_{\omega\to \infty}\dfrac{1}{\omega}=0$, following the principles established by Plotnikov \cite{P:1979}.

By employing the averaging of (\ref{eq:dotX_event}), we derive the following average system:
\begin{align}
\dfrac{dX_{\rm{av}}(\bar{t})}{d\bar{t}}&=\dfrac{1}{\omega}\mathcal{F}_{\rm{av}}\left(X_{\rm{av}}\right) \,, \label{eq:dotXav_event_1} \\
\mathcal{F}_{\rm{av}}\left(X_{\rm{av}}\right)&=\dfrac{1}{T}\int_{0}^{T}\mathcal{F}\left(\xi,X_{\rm{av}},0\right)d\xi \label{eq:mathcalFav_event} \,,
\end{align}
where we ``freeze'' the average states of $\hat{G}(\bar{t})$, $e(\bar{t})$, and $\tilde{\theta}(\bar{t})$, resulting in the following equations, for all $\bar t\in [\bar t_{\kappa}, \bar t_{\kappa+1})$:
\begin{align}
\frac{d\hat{G}_{\rm{av}}(\bar{t})}{d\bar{t}}&=\frac{1}{\omega}HK\hat{G}_{\rm{av}}(\bar{t})+\frac{1}{\omega}HKe_{\rm{av}}(\bar{t})\,, \label{eq:dotHatGav_event_1} \\
\frac{d\tilde{\theta}_{\rm{av}}(\bar{t})}{d\bar{t}}&=\frac{1}{\omega}KH\tilde{\theta}_{\rm{av}}(\bar{t})+\frac{1}{\omega}Ke_{\rm{av}}(\bar{t})\,, \label{eq:dotTildeThetaAv_event_1} \\
e_{\rm{av}}(\bar{t})&=\hat{G}_{\rm{av}}(\bar{t}_{\kappa})-\hat{G}_{\rm{av}}(\bar{t})\,, \label{eq:Eav_event_1} \\
\hat{G}_{\rm{av}}(\bar{t})&= H\tilde{\theta}_{\rm{av}}(\bar{t})\,, \label{eq:hatGav_event_1}
\end{align}
recalling that the matrix $HK$ is Hurwitz. Therefore, it is evident from (\ref{eq:dotHatGav_event_1}) that the ISS relationship of $\hat{G}_{\rm{av}}(\bar{t})$ with respect to the average measurement error $e_{\rm{av}}(\bar{t})$ holds. Thus, we can introduce the following ``Average NES Static-Triggering Condition'' for the average system.

\begin{definition}[\small{Average NES Static-Triggering Condition}] \label{def:averageStaticEvent}The Average Nash event-triggered condition consists of two components:
\begin{enumerate}
	\item The set of increasing sequences of time $I=\{I_{1}\,,I_{2}\}$ such that $I_{i}=\{\bar{t}^{i}_{0}\,, \bar{t}^{i}_{1}\,, \bar{t}^{i}_{2}\,,\ldots\}$ with $\bar{t}^{i}_{0}=0$, for all $i\in \{1,\ldots,N\}$, generated under the following rules:
		\begin{itemize}
			\item If $\left\{\bar{t} \in\mathbb{R}^{+}: \bar{t}>\bar{t}^{i}_{\kappa}  \wedge  \sigma_{i}|\hat{G}_{i}^{\rm{av}}(\bar{t})| -|e_{i}^{\rm{av}}(\bar{t})| <  0 \right\} = \emptyset$, then the set of the times of the events is $I_{i}=\{\bar{t}^{i}_{0}\,, \bar{t}^{i}_{1}\,, \ldots, \bar{t}^{i}_{\kappa}\}$.
			\item If $\left\{\bar{t} \in\mathbb{R}^{+}: \bar{t}>\bar{t}^{i}_{\kappa} \wedge  \sigma_{i}|\hat{G}_{i}^{\rm{av}}(\bar{t})| -|e_{i}^{\rm{av}}(\bar{t})| < 0 \right\} \neq \emptyset$, the next event time is given by
				\begin{align}
					\!\!\!\!\!\!\!\!\!\!\!\!\!\!\!\!\!\!\!\!\! \bar{t}^{i}_{\kappa+1}&=\inf\left\{\bar{t} \in\mathbb{R}^{+}: \bar{t}>\bar{t}^{i}_{\kappa} ~ \wedge ~ \sigma_{i}|\hat{G}_{i}^{\rm{av}}(\bar{t})| -|e_{i}^{\rm{av}}(\bar{t})| < 0 \right\}\,, \label{eq:tk+1_event_av}
				\end{align}
				\!\!\! consisting of the static event-triggering mechanism.
		\end{itemize}
	\item The $i$-th feedback control action updated at the  triggering instants such that
		\begin{align}
			u_{i}^{\rm av}(\bar{t})=K_{i}\hat{G}^{\rm av}_{i}(\bar{t})+K_{i}e^{\rm av}_{i}(\bar{t}) \,, \label{eq:U_MD4}
		\end{align}
		for all $\bar{t} \in \lbrack \bar{t}_{\kappa}\,, \bar{t}_{\kappa+1}\phantom{(}\!\!)$, $k\in \mathbb{N}$.
\end{enumerate} 
\end{definition}

\section{Stability Analysis}

The next theorem guarantees the local asymptotic stability of the NES control system employing static event-triggered execution, as depicted in Fig.~\ref{fig:blockDiagram_v3}.

\begin{theorem} \label{thm:NETESC_2}
Consider the closed-loop average dynamics of the gradient estimate (\ref{eq:dotHatGav_event_1}), the average error vector \eqref{eq:Eav_event_1}, and the average distributed event-triggering mechanism in Definition \ref{def:averageStaticEvent}. For a sufficiently large $\omega>0$, defined in (\ref{eq:new_omega}), the equilibrium $\hat{G}_{\rm{av}}(t)\equiv0$ is locally exponentially stable, and $\tilde{\theta}_{\rm{av}}(t)$ exponentially converges to zero. Furthermore, for the non-average system (\ref{eq:dotTildeTheta_v1_20240319}), there exist constants $m$ and $\bar{M}_{\theta}$ 
such that
\begin{align}
\|\theta(t)-\theta^{\ast}\| &\leq \bar{M}_{\theta}\exp(-mt)+\mathcal{O}\left(a+\frac{1}{\omega}\right)\,, \label{eq:normTheta_thm2} 
\end{align}
where $a=\sqrt{a_{1}^{2}+a_{2}^{2}}$, with $a_1$, $a_2$ defined in \eqref{eq:Si}--\eqref{eq:Mi}, and the constants $m$ and $\bar{M}_{\theta}$ 
depending on the triggering parameters $\sigma_1$, $\sigma_2$, and the initial condition $\theta(0)$, respectively. Additionally, there exists a lower bound $\tau^{\ast}$ for the inter-execution interval $t_{\kappa+1}-t_{\kappa}$ for all $\kappa \in \mathbb{N}$, preventing the Zeno behavior.
\end{theorem}

\textit{Proof:} The proof of the theorem is split into two main sections: stability analysis and avoidance of Zeno behavior.

\begin{flushleft}
\textcolor{black}{\underline{\it A. Stability Analysis}}
\end{flushleft}

Let us examine the following candidate Lyapunov function for the average system  (\ref{eq:dotHatGav_event_1}):
\begin{align}
V_{\text{av}}(\bar{t})=\hat{G}^{T}_{\text{av}}(\bar{t})P\hat{G}_{\text{av}}(\bar{t}) \,, \quad P=P^T>0\,. \label{eq:Vav_event_pf2}
\end{align} 
Since $HK$ is Hurwitz, given $Q=Q^T>0$ there exist $P=P^T$ such that the Lyapunov equation is $K^{T}H^{T}P+PHK=-Q$ and, therefore, the time derivative of (\ref{eq:Vav_event_pf2}) is given by  
\begin{align}
\frac{dV_{\text{av}}(\bar{t})}{d\bar{t}}&=-\frac{1}{\omega}\hat{G}_{\text{av}}^{T}(\bar{t})Q\hat{G}_{\text{av}}(\bar{t})+\frac{1}{\omega}e_{\text{av}}^{T}(\bar{t})H^{T}K^{T}P\hat{G}_{\text{av}}(\bar{t})\nonumber \\
&\quad+\frac{1}{\omega}\hat{G}_{\text{av}}^{T}(\bar{t})PKHe_{\text{av}}(\bar{t})\,, \label{eq:dotVav_event_1_pf2}
\end{align}
whose upper bound can be expressed as
\begin{align}
\frac{dV_{\text{av}}(\bar{t})}{d\bar{t}}&\leq\mathbb{-}\frac{\lambda_{\min}(Q)}{\omega}\|\hat{G}_{\text{av}}(\bar{t})\|^{2}\mathbb{+}
\frac{ 2\|PKH\|}{\omega}\|e_{\text{av}}(\bar{t})\| \|\hat{G}_{\text{av}}(\bar{t})\|. \label{eq:dotVav_event_2_pf2}
\end{align}
In the proposed event-triggering mechanism, the update law is given by (\ref{eq:tk+1_event_av}). Thus, by design $|e^{\rm{av}}_{1}(\bar{t})|\leq \sigma_{1}|\hat{G}^{\rm{av}}_{1}(\bar{t})|$, $|e^{\rm{av}}_{2}(\bar{t})|\leq \sigma_{2}|\hat{G}^{\rm{av}}_{2}(\bar{t})|$ and $\|e_{\rm{av}}(\bar{t})\|\leq \bar{\sigma}\|\hat{G}_{\rm{av}}(\bar{t})\|$, where $\bar{\sigma} = \max\left\{\sigma_{1},\sigma_{2}\right\}$. Additionally, the Rayleigh-Ritz inequality \cite{K:2002} give us $\lambda_{\min}(P)\|\hat{G}_{\text{av}}(\bar{t})\|^{2}\leq V_{\text{av}}(\bar{t}) \leq \lambda_{\max}(P)\|\hat{G}_{\text{av}}(\bar{t})\|^{2}$. Thus, inequality (\ref{eq:dotVav_event_2_pf2}) is upper bounded by 
\begin{align}
\frac{dV_{\text{av}}(\bar{t})}{d\bar{t}}&\leq-\frac{1}{\omega}\frac{\lambda_{\min}(Q)}{\lambda_{\max}(P)}\left(1-\frac{ 2\|PKH\|\bar{\sigma}}{\lambda_{\min}(Q)}\right)V_{\text{av}}(\bar{t}) \,, \label{eq:dotVav_event_5_pf2}
\end{align}
with $V_{\text{av}}(\bar{t})>0$ and $\dfrac{dV_{\text{av}}(\bar{t})}{d\bar{t}}<0$, for all $\bar{\sigma}<\dfrac{\lambda_{\min}(Q)}{ 2\|PHK\|}$. For instance, if we choose $\bar{\sigma}=\frac{\lambda_{\min}(Q)}{ 2\|PHK\|}\hat{\sigma}$, where $\hat{\sigma} \in (0,1)$, and defining $\alpha = \frac{\lambda_{\min}(Q)}{\lambda_{\max}(P)}$, inequality (\ref{eq:dotVav_event_5_pf2}) simply becomes
\begin{align}
\frac{dV_{\text{av}}(\bar{t})}{d\bar{t}}&\leq-\frac{\alpha\left(1-\hat{\sigma}\right)}{\omega}V_{\text{av}}(\bar{t}) \,. \label{eq:dotVav_20240307_1}
\end{align}
Using the Comparison Principle \cite[Lemma]{K:2002}, the solution of 
\begin{align}
\frac{d\bar{V}_{\text{av}}(\bar{t})}{d\bar{t}}=-\frac{\alpha\left(1-\hat{\sigma}\right)}{\omega}\bar{V}_{\text{av}}(\bar{t})\,, \quad \bar{V}_{\text{av}}(\bar{t}_{\kappa})=V_{\text{av}}(\bar{t}_{\kappa})
\end{align}
provides an upper bound $\bar{V}_{\text{av}}(\bar{t})$ for $V_{\text{av}}(\bar{t})$ such that
\begin{align}
V_{\text{av}}(\bar{t})\leq \bar{V}_{\text{av}}(\bar{t}) \,, \quad \forall \bar{t}\in \lbrack \bar{t}_{\kappa},\bar{t}_{\kappa+1}\phantom{(}\!\!) \,, \label{eq:VavBarVav_1_pf2}
\end{align}
with
\begin{align}
\bar{V}_{\text{av}}(\bar{t})=\exp\left(-\frac{\alpha\left(1-\hat{\sigma}\right)}{\omega}\bar{t}\right)V_{\text{av}}(\bar{t}_{\kappa})\,, \quad \forall \bar{t}\in \lbrack \bar{t}_{\kappa},\bar{t}_{\kappa+1}\phantom{(}\!\!)\,. \label{eq:_pf2}
\end{align}
By defining $\bar{t}_{\kappa}^{+}$ and $\bar{t}_{\kappa}^{-}$ as the right and left limits of $\bar{t}=\bar{t}_{\kappa}$, respectively, it is easy to verify that $V_{\text{av}}(\bar{t}_{\kappa+1}^{-})\leq \exp\left(-\dfrac{\alpha\left(1-\hat{\sigma}\right)}{\omega}(\bar{t}_{\kappa+1}^{-}-\bar{t}_{\kappa}^{+})\right)V_{\text{av}}(\bar{t}_{\kappa}^{+})$. Since $V_{\text{av}}(\bar{t})$ is continuous, $V_{\text{av}}(\bar{t}_{\kappa+1}^{-})=V_{\text{av}}(\bar{t}_{\kappa+1})$, $V_{\text{av}}(\bar{t}_{\kappa}^{+})=V_{\text{av}}(\bar{t}_{\kappa})$, and, consequently,
\begin{align}
    V_{\text{av}}(\bar{t}_{\kappa+1})\leq \exp\left(-\frac{\alpha\left(1-\hat{\sigma}\right)}{\omega}(\bar{t}_{\kappa+1}-\bar{t}_{\kappa})\right)V_{\text{av}}(\bar{t}_{\kappa})\,. \label{METES_eq:mmd_1_s}
\end{align}
Hence, for any $\bar{t}\geq 0$ in $ \bar{t}\in \lbrack \bar{t}_{\kappa},\bar{t}_{\kappa+1}\phantom{(}\!\!)$, $k \in \mathbb{N}$, one has 
\begin{small}
\begin{align}
   & V_{\text{av}}(\bar{t})\leq \exp\left(\mathbb{-}\frac{\alpha\left(1\mathbb{-}\hat{\sigma}\right)}{\omega}(\bar{t}\mathbb{-}\bar{t}_{\kappa})\right) V_{\text{av}}(\bar{t}_{\kappa}) \nonumber \\
    &\leq \exp\left(-\frac{\alpha\left(1\mathbb{-}\hat{\sigma}\right)}{\omega}(\bar{t}\mathbb{-}\bar{t}_{\kappa})\right) \exp\left(\mathbb{-}\frac{\alpha\left(1\mathbb{-}\hat{\sigma}\right)}{\omega}(\bar{t}_{\kappa}\mathbb{-}\bar{t}_{\kappa-1})\right)V_{\text{av}}(\bar{t}_{\kappa-1}) \nonumber \\
    &\leq \ldots \leq \nonumber \\
    &\leq \exp\left(\!\!\mathbb{-}\frac{\alpha\left(1\mathbb{-}\hat{\sigma}\right)}{\omega}(\bar{t}\mathbb{-}\bar{t}_{\kappa})\right) \prod_{i=1}^{i=k}\exp\left(\!\!\mathbb{-}\frac{\alpha\left(1\mathbb{-}\hat{\sigma}\right)}{\omega}(\bar{t}_{i}\mathbb{-}\bar{t}_{i-1})\right)V_{\text{av}}(\bar{t}_{i-1}) \nonumber \\
    &=\exp\left(\mathbb{-}\frac{\alpha\left(1\mathbb{-}\hat{\sigma}\right)}{\omega}\bar{t}\right) V_{\text{av}}(0)\,, \quad \forall \bar{t}\geq 0\,. \label{METES_eq:VavBarVav_2_pf2}
\end{align}
\end{small}
Now, by lower bounding the left-hand side and upper bounding the right-hand side of (\ref{METES_eq:VavBarVav_2_pf2}) with their counterparts in Rayleigh-Ritz inequality, we obtain:
\begin{align}
\lambda_{\min}(P)\|\hat{G}_{\text{av}}(\bar{t})\|^{2}&\leq \exp\left(-\frac{\alpha\left(1-\hat{\sigma}\right)}{\omega}\bar{t}\right)\lambda_{\max}(P)\|\hat{G}_{\text{av}}(0)\|^{2} \,. \label{eq:VavBarVav_3_pf2}
\end{align}
Then,
\begin{align}
\|\hat{G}_{\text{av}}(\bar{t})\|^{2}&\leq 
\left[\exp\left(-\frac{\alpha\left(1-\hat{\sigma}\right)}{2\omega}\bar{t}\right)\sqrt{\frac{\lambda_{\max}(P)}{\lambda_{\min}(P)}}\|\hat{G}_{\text{av}}(0)\|\right]^{2}\,, \label{eq:VavBarVav_4_pf2}
\end{align}
and
\begin{align}
\|\hat{G}_{\text{av}}(\bar{t})\|\leq\exp\left(-\frac{\alpha\left(1-\hat{\sigma}\right)}{2\omega}\bar{t}\right)\sqrt{\frac{\lambda_{\max}(P)}{\lambda_{\min}(P)}}\|\hat{G}_{\text{av}}(0)\|\,. \label{eq:normHatGav_1_pf2}
\end{align}
Once $H$ is invertible and, from (\ref{eq:hatGav_event_1}), $\tilde{\theta}_{\rm{av}}(\bar{t})=H^{-1}\hat{G}_{\rm{av}}(\bar{t})$, the following inequality is established 
\begin{align}
\|\tilde{\theta}_{\rm{av}}(\bar{t})\|&\leq\exp\left(-\frac{\alpha\left(1-\hat{\sigma}\right)}{2\omega}\bar{t}\right)M_{\theta}\|\tilde{\theta}_{\rm{av}}(0)\|\,, \label{eq:nomrTildeTheta_pf2} \\
M_{\theta}&=\sqrt{\frac{\lambda_{\max}(P)}{\lambda_{\min}(P)}}\|H^{-1}\| \|H\| \,.
\end{align}
Since (\ref{eq:dotTildeTheta_v1_20240319}) has a discontinuous right-hand side, but it is also $T$-periodic in $t$, and noting that the average system with state $\tilde{\theta}_{\text{av}}(\bar{t})$ is asymptotically stable according to (\ref{eq:nomrTildeTheta_pf2}), we can invoke the averaging theorem in \cite[Theorem~2]{P:1979} to conclude that
\begin{align}
\|\tilde{\theta}(t)-\tilde{\theta}_{\text{av}}(t)\|\leq\mathcal{O}\left(\frac{1}{\omega}\right)\,.
\end{align}
By applying the triangle inequality \cite{A:1957}, we also obtain:
\begin{align}
\|\tilde{\theta}(t)\|&\leq\|\tilde{\theta}_{\text{av}}(t)\|+\mathcal{O}\left(\frac{1}{\omega}\right)\nonumber \\
&\leq \exp\left(-\frac{\alpha\left(1-\hat{\sigma}\right)}{2}t\right)M_{\theta}\|\tilde{\theta}_{\text{av}}(0)\|+\mathcal{O}\left(\frac{1}{\omega}\right)\!. \label{eq:nomrTildeTheta_pf2_gold}
\end{align}
Now, from (\ref{eq:Si}) and Fig.~\ref{fig:blockDiagram_v3}, we can verify that
\begin{align}
\theta(t)-\theta^{\ast}=\tilde{\theta}(t)+S(t)\,, \label{eq:theta_3_event_pf2}
\end{align}
where $S(t):=[S_{1}(t)\,,S_{2}(t)]^{T}$ and whose the Euclidean norm satisfies
\begin{align}
&\|\theta(t)-\theta^{\ast}\|=\|\tilde{\theta}(t)+S(t)\| \leq \|\tilde{\theta}(t)\|+\|S(t)\| \nonumber \\
&\leq \exp\left(-\frac{\alpha\left(1-\hat{\sigma}\right)}{2}t\right)M_{\theta}\|\theta(0) - \theta^{\ast}\|+\mathcal{O}\left(a+\frac{1}{\omega}\right)\,, \label{eq:theta_4_event_pf2}
\end{align}
leading to inequality (\ref{eq:normTheta_thm2}), for appropriate $m$ and $\bar{M}_\theta$.

\begin{flushleft}
\textcolor{black}{\underline{\it B. Avoidance of Zeno Behavior}}
\end{flushleft}

Since the average closed-loop system consists of (\ref{eq:dotHatGav_event_1}), with the event-triggering mechanism  (\ref{eq:tk+1_event_av}) and the average control law (\ref{eq:U_MD4}), we can conclude that $\|e_{\text{av}}(\bar{t})\| \leq \bar{\sigma}\|\hat{G}_{\text{av}}(\bar{t})\|$, resulting in  
\begin{align}
\bar{\sigma}\|\hat{G}_{\text{av}}(\bar{t})\|^{2}-\|e_{\text{av}}(\bar{t})\|\|\hat{G}_{\text{av}}(\bar{t})\|\geq 0\,. \label{ineq:interEvents_1_static}
\end{align}
By using the Peter-Paul inequality \cite{W:1971}, $cd\leq \frac{c^2}{2\epsilon}+\frac{\epsilon d^2}{2}$ for all $c,d,\epsilon>0$, with $c=\|e_{\rm{av}}(\bar{t})\|$, $d=\|\hat{G}_{\rm{av}}(\bar{t})\|$ and $\epsilon=\bar{\sigma}$, the inequality (\ref{ineq:interEvents_1_static}) is lower bounded by
\begin{align}
&\bar{\sigma} \|\hat{G}_{\text{av}}(\bar{t})\|^{2}-\|e_{\text{av}}(\bar{t})\|\|\hat{G}_{\text{av}}(\bar{t})\|\geq \nonumber \\
&q\|\hat{G}_{\rm{av}}(\bar{t})\|^{2}-p\|e_{\rm{av}}(\bar{t})\|^2\,,\label{ineq:interEvents_2_static_pf2}
\end{align}
where  $q=\frac{\bar{\sigma}}{2}$ and $p=\frac{1}{2\bar{\sigma}}$. 
In \cite{G:2014}, it is shown that a lower bound for the inter-execution interval is given by the time duration it takes for the function
\begin{align}
\phi(\bar{t})=\sqrt{\frac{p}{q}}\frac{\|e_{\rm{av}}(\bar{t})\|}{\|\hat{G}_{\rm{av}}(\bar{t})\|} \label{eq:phi_1_static_pf2}
\end{align}
to go from 0 to 1. The time-derivative of (\ref{eq:phi_1_static_pf2}) is 
\begin{align}
\frac{d\phi(\bar{t})}{d\bar{t}}&=\sqrt{\frac{p}{q}}\frac{1}{\|e_{\rm{av}}(\bar{t})\|\|\hat{G}_{\rm{av}}(\bar{t})\|}\left[e_{\rm{av}}^{T}(\bar{t})\frac{de_{\rm{av}}(\bar{t})}{d\bar{t}}\right. \nonumber \\
&\quad\left.-\hat{G}_{\rm{av}}^{T}(\bar{t})\frac{d\hat{G}_{\rm{av}}(\bar{t})}{d\bar{t}}\left(\frac{\|e_{\rm{av}}(\bar{t})\|}{\|\hat{G}_{\rm{av}}(\bar{t})\|}\right)^2\right]\,. \label{eq:dotPhi_1_static_pf2}
\end{align}
Then, the following estimate holds:
\begin{align}
\frac{d\phi(\bar{t})}{d\bar{t}}&\leq\frac{\|HK\|}{\omega}\sqrt{\frac{p}{q}}\left\{1+2\frac{\|e_{\rm{av}}(\bar{t})\|}{\|\hat{G}_{\rm{av}}(\bar{t})\|}+\frac{\|e_{\rm{av}}(\bar{t})\|^2}{\|\hat{G}_{\rm{av}}(\bar{t})\|^2}\right\}\,. \label{eq:dotPhi_2_static_pf2}
\end{align}

Hence,  using (\ref{eq:phi_1_static_pf2}), inequality (\ref{eq:dotPhi_2_static_pf2}) is rewritten as
\begin{align}
\omega\frac{d\phi(\bar{t})}{d\bar{t}}&\leq\|HK\|\sqrt{\frac{p}{q}}+2\|HK\|\phi(\bar{t})+\|HK\|\sqrt{\frac{q}{p}}\phi^{2}(\bar{t})\,. \label{eq:dotPhi_3_static_pf2}
\end{align}
From  the time-scaling $t =\frac{\bar{t}}{\omega}$, inequality (\ref{eq:dotPhi_2_static_pf2}) and invoking the Comparison Lemma \cite{K:2002}, a lower bound for the inter-execution time is found as
\begin{align}
\tau^{\ast}=\int_{0}^{1}\dfrac{1}{b_{0}+b_{1}\xi+b_{2}\xi^2}{\it d}\xi\,,\label{eq:tauAst_static_pf2}
\end{align}
with $b_{0}=\dfrac{\|HK\| }{\bar{\sigma}}$, $b_{1}=2\|HK\|$ and $b_{2}=\bar{\sigma}\|HK\|$. Therefore, the Zeno behavior is avoided not only for the average closed-loop system in time $\bar{t}$ but also for the original system in time $t$ since $\bar{t}={\omega}t$ is only a time compression (dilation) for $\omega$ sufficiently large (small), which means that $\bar{\tau}^{\ast}\mathbb{=}{\omega}\tau^{\ast}$ will still be a finite number, establishing a minimum switching time to rule out any 
Zeno behavior in time $t\in\mathbb{R}_{+}$ or $\bar{t}\in\mathbb{R}_{+}$. \hfill $\square$

\section{SIMULATION RESULTS}

This section shows simulation results in order to illustrate the distributed NES based on event-triggering mechanism. The investigated system captures a noncooperative game involving two firms operating in an duopoly market framework. These firms engage in competition aimed at maximizing their profits $J_{i}(t)$ through the pricing strategy $u_{i}(t)$ of their respective products, without sharing any information between the players. The plant parameters are consistent with those outlined in \cite{FKB:2012}: initial conditions are set as $\hat{\theta}_{1}(0)=50$ and $\hat{\theta}_{2}(0)=110/3$. Additionally, $S_{d} = 100$, $p=0.20$, $m_{1}=30$ and $m_{2}=30$. According to (\ref{eq:NE_v2}), these parameters lead us to a unique Nash equilibrium given by:
\begin{align}
\theta^{\ast}&=[43.3333\,, ~36.6667]^{T}\,,  \label{eq:theta_NE} \\
J^{\ast}&=[888.8889\,, ~222.2222]^{T}\,. \label{eq:J_NE}
\end{align}
The controller parameters are chosen as follows: for the NES strategy, we set $a_{1}=0.075$, $a_{2}=0.050$, $K_{1}=2$, $K_{2}=5$, $\omega_{1}=27$, and $\omega_{2}=22$. Regarding the distributed event-triggering mechanism, the parameters are $\sigma_{1}=0.85$ and $\sigma_{2}=0.95$.

\begin{figure}[h!]
	\centering
	\subfigure[Aperiodic sample-and-hold control inputs. \label{fig:u}]{\includegraphics[width=4.15cm]{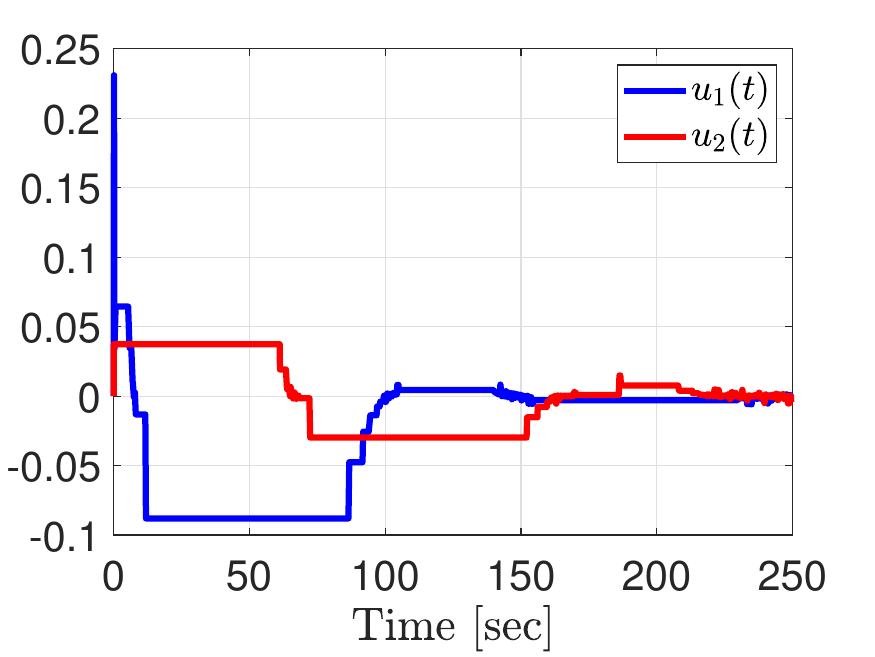}}
	~
	\subfigure[Time evolution of the controllers updates. \label{fig:controlUpdate}]{\includegraphics[width=4.15cm]{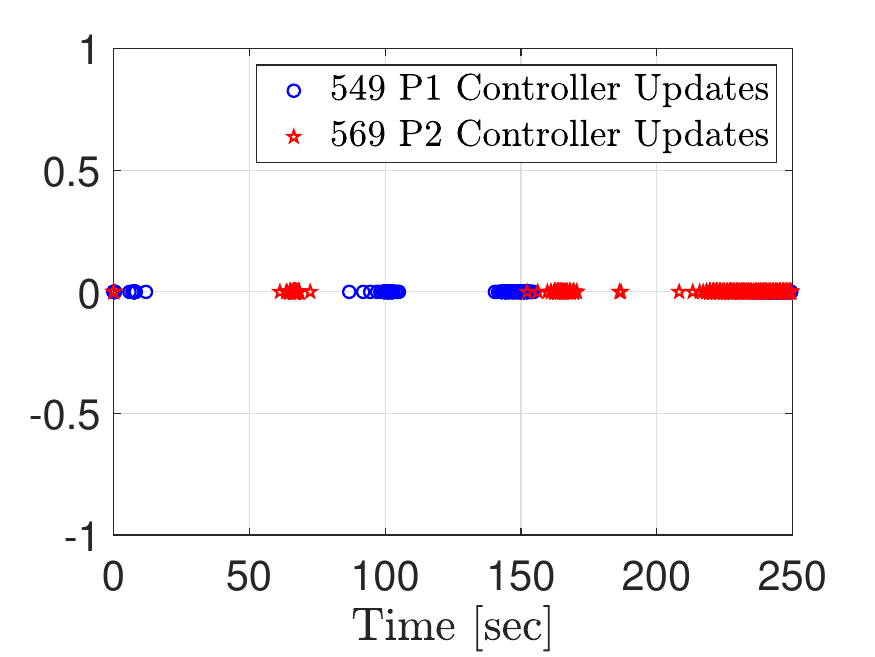}}
	\\
	\subfigure[Input signals of payoff functions, $\theta(t)$. \label{fig:theta}]{\includegraphics[width=4.15cm]{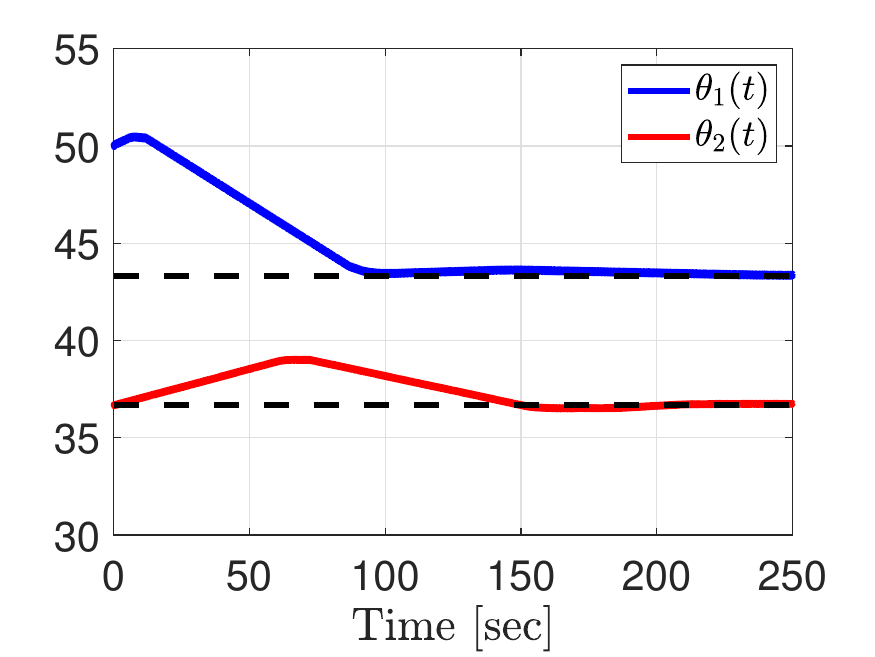}}
	~
	\subfigure[Payoff functions. \label{fig:J}]{\includegraphics[width=4.15cm]{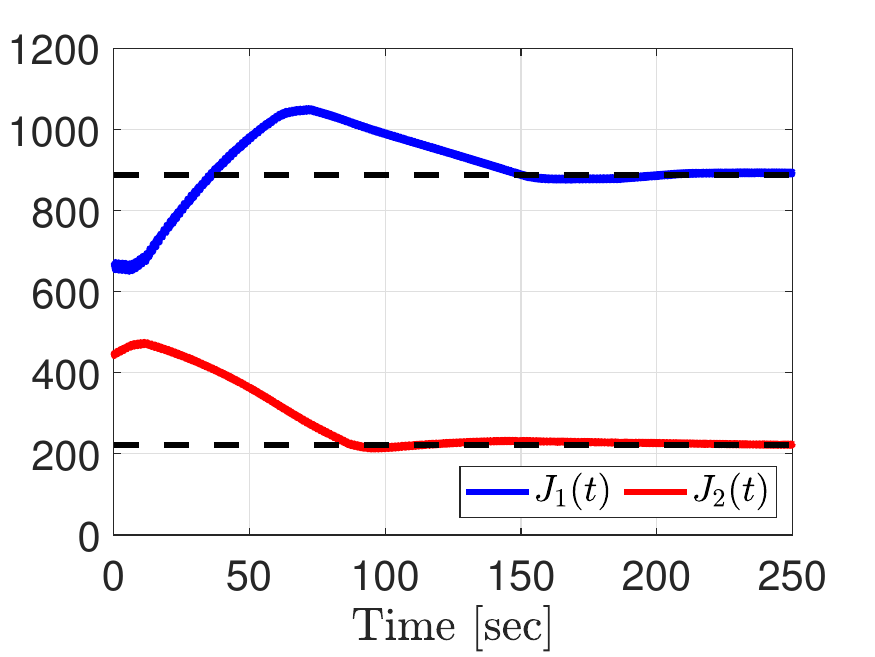}}
	\caption{Numerical simulations for NES in noncooperative duopoly games through distributed ETC policies. }\label{fig:ET_NES}
\end{figure}

To achieve the Nash equilibrium (\ref{eq:theta_NE}) with optimal payoffs (\ref{eq:J_NE}), players P1 and P2 implement the proposed decentralized event-triggered NES strategy employing sinusoidal perturbations to dictate their optimal actions. The time evolution of the proposed approach is depicted in Fig. \ref{fig:ET_NES}. Fig. \ref{fig:u} shows the aperiodic update behavior of the players' actions. Recalling that each player estimates a distinct gradient component and decides independently when to trigger and update the corresponding player action, these updates occur autonomously. In a simulation spanning 250 seconds, the actions of player P1 were updated 549 times and those of player P2 updated 569 times---see Fig.~\ref{fig:controlUpdate}. Notably, while these updates are independent, they collectively drive the system towards the Nash equilibrium, as depicted in Fig. \ref{fig:theta}. Moreover, Fig. \ref{fig:J} illustrates how, within the described duopoly market structure without information sharing, each player can maximize their profits $J_{i}(t)$ by employing the proposed decentralized NES strategy via ETC to set the product prices $\theta_{1}(t)$ and $\theta_{2}(t)$ in the noncooperative game.

\section{CONCLUSIONS} 

This paper introduced a novel method for achieving locally stable convergence to Nash equilibrium in noncooperative games through a distributed event-triggered control scheme using extremum seeking. This strategy enables the players of a duopoly game to improve their profits without sharing information, with action updates occurring independently since each player estimates and updates her actions based on different gradient components. The proposed approach contributes not only to game theory but also offers a practical method for decentralized decision-making in complex systems, demonstrating effectiveness even under limited bandwidth and preserving closed-loop stability.

\begin{small}
 
\end{small}

\end{document}